\documentclass[12 pt]{amsart}

\usepackage{amsmath,amsthm,amssymb,shuffle,stmaryrd,mathtools}
\usepackage{subcaption}
\DeclareCaptionSubType * [Alph]{table}
\usepackage{ytableau}
\usepackage{enumitem}
\usepackage[hmargin=1 in, vmargin = 1 in]{geometry}
\usepackage{hyperref}
\usepackage{tikz-cd}
\tikzcdset{every label/.append style = {font = \footnotesize}}

\usepackage[protrusion=true]{microtype}

\definecolor{darkblue}{rgb}{0.0,0,0.7}

\newcommand{\newword}[1]{\textcolor{darkblue}{\textbf{\emph{#1}}}}

\newcommand{\Z}{\mathbb{Z}}

\newcommand{\Plactic}{{\mathbf{P}}}
\newcommand{\Free}{\mathbf{F}}
\newcommand{\sPlactic}{\mathbf{S}}
\newcommand{\Monoid}{\mathbf{M}}
\newcommand{\Abelian}{\mathbf{C}}

\newcommand{\Category}{\underline{\mathsf{SPlac}}}
\newcommand{\LSCategory}{\underline{\mathsf{Plac}}}

\newcommand{\PlacticSchur}{\mathcal{S}}
\newcommand{\FreeSchur}{\widehat{\mathcal{S}}}

\newcommand{\sFreeSchur}{\widehat{\mathcal{P}}}

\newcommand{\symmetricgroup}{\mathfrak{S}}

\newcommand{\reading}{\mathsf{rw}}

\newcommand{\cA}{\mathcal{A}}
\newcommand{\cB}{\mathcal{B}}
\newcommand{\cC}{\mathcal{C}}
\newcommand{\cN}{\mathcal{N}}

\newcommand{\SSYT}{\mathrm{SSYT}}
\newcommand{\SYT}{\mathrm{SYT}}
\newcommand{\shSSYT}{\mathrm{ShSSYT}}

\newcommand{\id}{\mathrm{id}}

\newcommand{\hookword}{\mathsf{hook}}

\newtheorem{thm}{Theorem}[section]
\newtheorem{theorem}[thm]{Theorem}

\newtheorem{proposition}[thm]{Proposition}

\theoremstyle{remark}
\newtheorem{definition}[thm]{Definition}
\newtheorem{remark}[thm]{Remark}
\newtheorem{examplex}[thm]{Example}
\newenvironment{example}
  {\pushQED{\qed}\examplex}
  {\popQED\endexamplex}

\numberwithin{equation}{section}

\usepackage[colorinlistoftodos]{todonotes}

\AtBeginDocument{%
   \def\MR#1{}
}

\begin{document}

\title{A universal characterization of the shifted plactic monoid}

\author{Santiago Estupi\~n\'an-Salamanca}
\address[SE]{Department of Combinatorics \& Optimization, University of Waterloo, Waterloo ON, Canada}
\email{sestupinan@uwaterloo.ca}
\author{Oliver Pechenik}
\address[OP]{Department of Combinatorics \& Optimization, University of Waterloo, Waterloo ON, Canada}
\email{opecheni@uwaterloo.ca}

\date{\today}


\begin{abstract}
The \emph{plactic monoid} $\Plactic$ of Lascoux and Sch\"{u}tzenberger (1981) plays an important role in proofs of the Littlewood--Richardson rule for computing multiplicities in the linear representation theory of the symmetric group $\symmetricgroup_n$ and the cohomology of Grassmannians. Commonly, $\Plactic$ is defined as a quotient of a free monoid by relations derived from a careful analysis of Schensted's insertion algorithm and the jeu de taquin algorithm on semistandard Young tableaux. However, Lascoux and Sch\"{u}tzenberger also gave an intrinsic characterization of $\Plactic$ via a universal property.

Serrano's (2010) \emph{shifted plactic monoid} $\sPlactic$ is an analogue of $\Plactic$ that governs instead the projective representation theory of $\symmetricgroup_n$ and the cohomology of isotropic Grassmannians. We provide a universal property for $\sPlactic$, analogous to the Lascoux--Sch\"{u}tzenberger characterization of $\Plactic$.
\end{abstract}

\maketitle

\section{Introduction}

The \emph{plactic monoid}, first discovered by D.~Knuth \cite{Knuth} and extensively studied by A.~Lascoux--M.-P.~Sch\"utzenberger \cite{Lascoux.Schutzenberger:plaxique} (see also, \cite[\S 5]{Lothaire}), plays a fundamental role in symmetric function theory and its applications both to the representation theory of symmetric and general linear groups and to the Schubert calculus of Grassmannians. The plactic monoid is a quotient of a free monoid by a pair of homogeneous degree $3$ equivalences (\emph{Knuth relations}), extracted from Knuth's \cite{Knuth} careful combinatorial analysis of insertion algorithms for Young tableaux. 

In \cite[Th\'eor\`eme~2.16]{Lascoux.Schutzenberger:plaxique}, Lascoux--Sch\"utzenberger provide an intrinsic characterization of the plactic monoid as the largest quotient satisfying a few reasonable properties of an essentially algebraic nature. From this perspective, the Knuth relations appear spontaneously, in a fashion we find enlightening. In particular, this approach essentially avoids tableau combinatorics. The first result (Theorem~\ref{thm:LS}) of this paper is a rephrasing of \cite[Th\'eor\`eme~2.16]{Lascoux.Schutzenberger:plaxique} in categorical language, identifying the plactic monoid as the initial object of an appropriate category that we construct. Although this result is merely a repackaging of \cite[Th\'eor\`eme~2.16]{Lascoux.Schutzenberger:plaxique}, we provide a complete proof. One motivation for completeness is that the original proof from \cite{Lascoux.Schutzenberger:plaxique} is incomplete and contains several inaccuracies. Our second motivation is that the argument  serves as a warm-up for the proof of our main result (Theorem~\ref{thm:main}), which we now turn to describing.


Recall that the plactic monoid plays a key role in the linear representation theory of symmetric groups and the Schubert calculus of Grassmannians. Of similar importance for the \emph{projective} representation theory of symmetric groups and the Schubert calculus of \emph{isotropic} Grassmannians is the \emph{shifted plactic monoid} of L.~Serrano \cite{Serrano}; for background on projective representation theory and the associated combinatorics, see \cite{Hoffman.Humphreys}, while for the geometric perspective we recommend \cite[\S 6]{Gillespie}. The shifted plactic monoid, like the ordinary plactic, is a quotient of a free monoid by certain homogeneous equivalences, in this case $8$ relations of degree $4$. Serrano obtained these relations through a somewhat complicated analysis of the relationship between jeu de taquin for shifted tableaux \cite{Worley, Sagan} and M.~Haiman's \emph{mixed insertion} \cite{Haiman}. As evidence of the non-obviousness of this construction, we point to the over $30$ year gap between development of the necessary combinatorial tools by Sagan, Worley, and Haiman and the invention of the shifted plactic monoid by Serrano.

Our main result is a characterization of the shifted plactic monoid by universal property. This characterization is analogous to the Lascoux--Sch\"utzenberger characterization of the ordinary plactic monoid. (An similarly analogous characterization of the \emph{hypoplactic monoid} was given in \cite[Theorem~6.13]{Novelli}.)  Moreover, our approach avoids consideration of shifted jeu de taquin or mixed insertion. 

\medskip
We now give precise statements of these results.
An \newword{alphabet} $\cA$ is a totally ordered set. A \newword{morphism} of alphabets is a function $f : \cA \to \cB$ such that if $\alpha_1, \alpha_2 \in \cA$ with $\alpha_1 < \alpha_2$, then $f(\alpha_1) < f(\alpha_2)$. Let $\Free(\cA)$ denote the free monoid on the alphabet $\cA$. A morphism $f : \cA \to \cB$ of alphabets induces an \newword{ordered morphism} $\Free(f) : \Free(\cA) \to \Free(\cB)$. Let $\Abelian(\cA)$ denote the free commutative monoid on the alphabet $\cA$ and let $\epsilon : \Free(\cA) \to \Abelian(\cA)$ be the natural abelianization quotient map.

\begin{remark}\label{rem:completion}
All the monoids we consider turn out to be graded. For a graded monoid $\Monoid = \Free(\cA)/ \sim$ and $\cB$ a finite subalphabet of $\cA$, we write $\Monoid(\cB)$ for the image of $\Free(\cB)$ under the quotient map and we abuse notation by writing $\Z \Monoid$ for the graded inverse limit of the monoid algebras $\Z\Monoid(\cB)$ for finite subalphabets $\cB$. In particular, if $\cA$ is infinite, $\Z \Monoid$ generally allows consideration of certain infinite $\Z$-linear combinations, while the monoid algebra of $\Monoid$ certainly does not.
\end{remark}

We say $I \subseteq \cA$ is an \newword{interval} if $i < k < j \in \cA$ with $i,j \in I$ implies that $k \in I$. For $I$ an interval of $\cA$, let $\rho_I : \Free(\cA) \to \Free(I)$ be the map that deletes all letters outside of the interval $I$. 
We extend $\rho_I$ linearly to a map $\rho_I : \Z \Free(\cA) \to \Z \Free(I)$ between the corresponding monoid algebras.
Let $\kappa : \Free(\cA) \to \Plactic(\cA)$ be the projection map given by quotienting by the Knuth relations (see Definition~\ref{def:Knuth}). Let $\lambda : \Plactic(\cA) \to \Abelian(\cA)$ be the abelianization map.

Consider the category $\LSCategory(\cA)$ of monoids $\Monoid$ equipped with morphisms
\[\phi : \Free(\cA) \twoheadrightarrow \Monoid \quad \text{and} \quad \psi : \Monoid \to \Abelian(\cA)
\]such that 
\begin{enumerate}[label=\textrm{(Plac.\arabic*)}]
    \item $\psi \circ \phi = \epsilon$; \label{plac1}
    \item \ytableausetup{boxsize=0.3em} the images $\phi(\FreeSchur_{\ydiagram{1}}), \phi(\FreeSchur_{\ydiagram{1,1}}) \in \Z\Monoid$ of the \emph{free Schur functions} (see Section~\ref{sec:free!}) \[
    \{ \FreeSchur_{\ydiagram{1}}, \FreeSchur_{\ydiagram{1,1}} \} \subset \Z\Free(\cA)
    \]
    commute with each other; \label{plac2}
    \item for any finite subalphabets $\cB, \cC$ of $\cA$, any ordered morphism $\omega : \Free(\cB) \to \Free(\cC)$, and any $w_1, w_2 \in \Free(\cB)$, if $\phi(w_1) = \phi(w_2)$, then $(\phi \circ \omega)(w_1) = (\phi \circ \omega)(w_2)$; and \label{plac3}
    \item for any interval $I \subseteq \cA$ and any $w_1, w_2 \in \Free(\cA)$, if $\phi(w_1) = \phi(w_2)$, then \[(\phi \circ \rho_I)(w_1) = (\phi \circ \rho_I)(w_2).\] \label{plac4}
\end{enumerate}
A morphism $\theta : (\Monoid', \phi', \psi') \to (\Monoid, \phi, \psi)$ in $\LSCategory(\cA)$ is a morphism of monoids such that the diagram
\begin{equation}
\label{eq:LS_rhombus}
\begin{tikzcd}[column sep=huge]
& \Monoid' \ar[dr, "\psi'"] \ar[dd,"\theta"'] 
& \\
\Free(\cA) \ar[ur, "\phi'"] \ar[dr, "\phi"']
&
& \Abelian(\cA) \\
& \Monoid \ar[ur, "\psi"']
&
\end{tikzcd}
\end{equation}
commutes.

\begin{theorem}\label{thm:LS}
The plactic monoid $(\Plactic(\cA), \kappa, \lambda)$ is the initial object of the category $\LSCategory(\cA)$.
\end{theorem}

In addition to being phrased in less categorical language, the theorem of Lascoux--Sch\"utzenberger \cite[Th\'eor\`eme~2.16]{Lascoux.Schutzenberger:plaxique} replaces \ref{plac2} with commutation of free Schur functions for \emph{all} column shapes; however, inspection of their proof reveals that only the case stated in \ref{plac2} is necessary.
We prove Theorem~\ref{thm:LS} in Section~\ref{sec:LS_universal_proof}. The main content of the proof is showing, without reference to an insertion algorithm, that the Knuth relations hold in every $\Monoid(\cA) \in \LSCategory$; this is done primarily by analyzing the commutativity condition \ref{plac2}. In particular, the theorem and its proof only require knowing the definitions of free Schur functions in the very smallest cases $\FreeSchur_{\ydiagram{1}}$ and $\FreeSchur_{\ydiagram{1,1}}$, which are natural to define even in the absence of a full tableau theory.
Note that the final object of $\LSCategory(\cA)$ is easily seen to be $(\Abelian(\cA),\epsilon, \id)$.

\bigskip
Consider now the category $\Category(\cA)$ of monoids $\Monoid$ equipped with morphisms 
\[\phi : \Free(\cA) \twoheadrightarrow \Monoid \quad \text{and} \quad \psi : \Monoid \to \Plactic(\cA)
\]such that 
\begin{enumerate}[label=\textrm{(SPlac.\arabic*)}]
    \item $\psi \circ \phi = \kappa$; \label{splac1}
    \item \ytableausetup{boxsize=0.3em} the images $\phi(\sFreeSchur_{\ydiagram{1}}), \phi(\sFreeSchur_{\ydiagram{2,1+1}}) \in \Z\Monoid$ of the \emph{shifted free Schur functions} (see Section~\ref{sec:free!})
    \[
    \{ \sFreeSchur_{\ydiagram{1}}, \sFreeSchur_{\ydiagram{2,1+1}} \} \subset \Z\Free(\cA)
    \]
    commute; \label{splac2}
    \item for any finite subalphabets $\cB, \cC$ of $\cA$, any ordered morphism $\omega : \Free(\cB) \to \Free(\cC)$, and any $w_1, w_2 \in \Free(\cB)$, if $\phi(w_1) = \phi(w_2)$, then $(\phi \circ \omega)(w_1) = (\phi \circ \omega)(w_2)$; and \label{splac3}
    \item for any interval $I \subseteq \cA$ and any $w_1, w_2 \in \Free(\cA)$, if $\phi(w_1) = \phi(w_2)$, then \[(\kappa \circ \rho_I)(w_1) = (\kappa \circ \rho_I)(w_2).\] \label{splac4}
\end{enumerate}
Note that while objects of $\LSCategory(\cA)$ naturally sit between $\Free(\cA)$ and $\Abelian(\cA)$, objects of $\Category(\cA)$ naturally sit between $\Free(\cA)$ and $\Plactic(\cA)$. The defining properties of $\Category(\cA)$ are mostly given by directly adapting those for $\LSCategory(\cA)$ to this setting. The main surprise is the appearance of $\kappa$ in \ref{splac4}, where one might naturally expect $\phi$ by analogy.
A morphism $\theta : (\Monoid', \phi', \psi') \to (\Monoid, \phi, \psi)$ in $\Category(\cA)$ is a morphism of monoids such that the diagram
\begin{equation}
 \label{eq:shifted_rhombus}
\begin{tikzcd}[column sep=huge]
& \Monoid' \ar[dr, "\psi'"] \ar[dd,"\theta"'] 
& \\
\Free(\cA) \ar[ur, "\phi'"] \ar[dr, "\phi"']
&
& \Plactic(\cA) \\
& \Monoid \ar[ur, "\psi"']
&
\end{tikzcd}
\end{equation}
commutes.

Let $\sigma : \Free(\cA) \to \sPlactic(\cA)$ be the projection map given by quotienting by the \emph{shifted Knuth relations} (see Definition~\ref{def:shifted_plactic}) and let $\pi : \sPlactic(\cA) \to \Plactic(\cA)$ be the projection map given by quotienting by the ordinary Knuth relations. Comparing Definitions~\ref{def:Knuth} and~\ref{def:shifted_plactic} makes clear that $\kappa = \pi \circ \sigma$. 

Our main result is the following, describing how the shifted plactic monoid $\sPlactic(\cA)$ may be characterized without reference to the shifted Knuth relations. Our approach also avoids consideration of Haiman's \emph{mixed insertion} and is analogous to the Lascoux--Sch\"utzenberger characterization of $\Plactic(\cA)$. 
\begin{theorem}\label{thm:main}
The shifted plactic monoid $(\sPlactic(\cA), \sigma, \pi)$ is the initial object of the category $\Category(\cA)$.
\end{theorem}

We prove Theorem~\ref{thm:main} in Section~\ref{sec:shifted_universal_proof}.  
Note that the final object of is $\Category(\cA)$ is just $(\Plactic(\cA),\kappa, \id)$. We also remark that Theorem~\ref{thm:main} and its proof only require knowing the definitions of shifted free Schur functions in the very small cases $\FreeSchur_{\ydiagram{1}}$ and $\FreeSchur_{\ydiagram{2,1+1}}$, which are both natural to define from a hook word perspective, even in the absence of a tableau theory.

\bigskip
{\bf This paper is organized as follows.} In Section~\ref{sec:background}, we recall necessary background and motivation, and we give some new definitions. In Section~\ref{sec:LS_universal_proof}, we prove Theorem~\ref{thm:LS}, following the strategy of Lascoux--Sch\"utzenberger \cite{Lascoux.Schutzenberger:plaxique}. In Section~\ref{sec:shifted_universal_proof}, we prove the main result, Theorem~\ref{thm:main}. In Section~\ref{sec:staircase}, we explain how one obtains the same theories by axiomatizing the categories $\LSCategory(\cA)$ and $\Category(\cA)$ using, respectively,
 $\FreeSchur_{\ydiagram{2}}$ in place of $\FreeSchur_{\ydiagram{1,1}}$ and $\sFreeSchur_{\ydiagram{3}}$ in place of $\sFreeSchur_{\ydiagram{2,1+1}}$. 

\ytableausetup{boxsize=normal}

\section{Background and definitions}\label{sec:background}

\subsection{Young tableaux}
Fix an alphabet $\cA$.
A \newword{partition} $\nu = (\nu_1 \geq \nu_2 \geq \dots \geq \nu_k)$ is a finite weakly decreasing list of positive integers. The number $k$ is called the \newword{length} $\ell(\nu)$ of the partition $\nu$. We write $|\nu| = \sum_i \nu_i$. Given partitions $\nu, \mu$, we write $\mu \subseteq \nu$ if $\ell(\mu) \leq \ell(\nu)$ and we have $\mu_i \leq \nu_i$ for all $1 \leq i \leq \ell(\mu)$.

The \newword{Young diagram} of a partition $\nu$ is given by left-justifying $\nu_i$ boxes in row $i$, where we treat row $1$ as the top row. 
For example, 
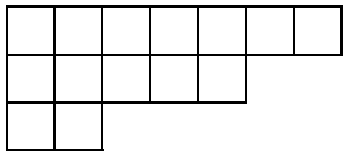
\[
\ydiagram{7,5,2}
\] 
is the Young diagram for the partition $\beta = (7,5,2)$. If $\mu \subseteq \nu$, $\nu / \mu$ denotes the set-theoretic difference of the Young diagrams.
A \newword{semistandard Young tableau} of \newword{shape} $\nu$ is a filling of the boxes of the Young diagram of $\nu$ by letters of $\cA$ such that the box labels weakly increase from left to right along rows and strictly increase from top to bottom down columns. For example,
\[
\ytableaushort{1112467,25555,49}
\]
is a semistandard Young tableau of shape $\beta$ in the alphabet $\cN \coloneqq \{1 < 2 < \dots \}$. We write $\SSYT(\nu)$ for the set of all semistandard tableaux of shape $\nu$. (The choice of alphabet should always be clear from context.) The \newword{content} of a tableau $T$ is the function $c_T : \cA \to \cN$ such that $c_T(a)$ equals the number of boxes filled with the letter $a$. A tableau $T \in \SSYT(\nu)$ in the alphabet $\cN$ is \newword{standard} if $c_T(i) = 1$ for $i \leq |\nu|$ and $c_T(i) = 0$ for $i > |\nu|$.
For example,
\[
\ytableaushort{1346789,
2{10}{11}{13}{14},
5{12}}
\]
is a standard tableau of shape $\beta$. We write $\SYT(\nu)$ for the set of standard tableaux of shape $\nu$. Returning to the arbitrary alphabet $\cA$, the \newword{Schur function} $s_\nu$ is the symmetric function 
\[
s_\nu \coloneqq \sum_{T \in \SSYT(\nu)} {\bf x}^{c_T} = \sum_{T \in \SSYT(\nu)} \prod_{a \in \cA} x_a^{c_T(a)} \in \Z\llbracket x_a : a \in \cA \rrbracket_{\rm gr},
\]
an element of the ring of power series of bounded degree.

A partition $\nu$ is \newword{strict} if $\nu_1 > \nu_2 > \dots > \nu_{\ell(\nu)}$. A strict partition has a \newword{shifted Young diagram} where the boxes of row $i$ are indented $i-1$ positions to the right. For example, 
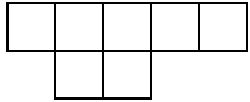
\[
\ydiagram{5,1+2}
\]
is the shifted Young diagram of the strict partition $\eta = (5,2)$. Consider the doubled alphabet 
$\cA' \coloneqq \bigsqcup_{a \in \cA} \{a, a'\}$ where for $a< b$, we set $a' < a < b'$.
A \newword{shifted semistandard (Young) tableau} of strict partition shape $\nu$ is a filling of the boxes of the shifted Young diagram for $\nu$ by elements of $\cA'$ such that 
\begin{itemize}
    \item the box labels weakly increase from left to right along rows and from top to bottom down columns,
    \item each $a'$ appears at most once in any row,
    \item each $a$ appears at most once in any column, and
    \item the leftmost box of each row contains an unprimed letter.
\end{itemize}
For example,
\[
\ytableaushort{11{3'}{4}{5'},\none 2{3'}}
\]
is a shifted semistandard Young tableau of shifted shape $\eta$ in the alphabet $\cN'$. We write $\shSSYT(\nu)$ for the set of all shifted semistandard Young tableaux of strict partition shape $\nu$. 
The \newword{content} of a shifted semistandard tableau $T$ in the doubled alphabet $\cA'$ is the function $c_T : \cA \to \cN$ where $c_T(a)$ denotes the number of boxes that are filled with either the letter $a'$ or $a$. 
The \newword{$P$-Schur function} $P_\nu$ is the symmetric function 
    $$P_\nu \coloneqq \sum_{T\in \shSSYT(\nu)}{\bf x}^{c_T}=\sum_{T\in \shSSYT(\nu)} \prod_{a \in \cA}x_a^{c_T(a)}.$$


\subsection{Free and plactic Schur functions}\label{sec:free!}
In \cite{Schutzenberger}, Sch\"utzenberger introduced a new approach to symmetric function theory, founded on a now prominent algorithm known as \emph{jeu de taquin}. This contribution, leading to the first complete proof of the Littlewood--Richardson rule, strengthened the connections between the combinatorial theory of Young tableaux and the structure of Schur functions. In particular, it showed that Schur functions can be lifted to the noncommutative algebras $\Z\Plactic(\cA)$ and $\Z \Free(\cA)$. Precisely, for $\nu$ a partition of length $\ell$, we define the \newword{free Schur function} $\FreeSchur_\nu$ as 
\begin{equation}\label{eq:FreeSchur}
    \FreeSchur_\nu \coloneqq \sum_{ w=w_1\ldots w_\ell } w \in \Z \Free(\cA),
\end{equation}
where the sum is over all concatenations of weakly increasing segments $w=w_1\ldots w_\ell$ such that $w_i$ is a longest weakly increasing subword of $w_{i-1}w_i$, and $w_i$ is of length $\nu_{\ell+1 - i}$.

\begin{remark}
    Our definition of free Schur functions is a variation of a special case of that from \cite{Lascoux.Schutzenberger:plaxique}, which is defined for finite alphabets and for particular choices of standard tableaux $T \in \SYT(\nu)$. In this special case, we can repackage the definition as above to remove any mention of tableaux. Nonetheless, we have 
    $\epsilon(\FreeSchur_T) = \epsilon(\FreeSchur_{\nu})$
    where $\FreeSchur_T$ stands for the Lascoux--Sch\"utzenberger free Schur functions. 
    Moreover, in the finite alphabet case, there exists a standard Young tableau $R$ so that 
    $\FreeSchur_\nu = \FreeSchur_R.$
    Note that for general $T$, the Lascoux--Sch\"{u}tzenberger free Schur function $\FreeSchur_T$ does not correspond to any of the free Schur functions indexed by a partition $\FreeSchur_\nu$. See also, \cite{Fomin.Greene} for similar ideas.
\end{remark}

It turns out that $\epsilon(\FreeSchur_\nu) = s_\nu$, where $\epsilon$ is the abelianization map. One of the chief virtues of this approach, developed in greater depth in \cite{Lascoux.Schutzenberger:plaxique,Schutzenberger:pourleplaxique,Lothaire}, is that it breaks apart the original summands of $s_\nu$, having the effect of making them more docile. For instance, this approach led to the first proof of a combinatorial description of the \newword{Littlewood--Richardson coefficients} $c_{\nu,\mu}^\xi$ defined by
$$s_\nu s_\mu=\sum_{\xi} c_{\nu,\mu}^\xi s_\xi.$$

An important tool in the proof of the Littlewood--Richardson rule is the following.
\begin{definition}[\cite{Knuth,Lascoux.Schutzenberger:plaxique}]\label{def:Knuth}
    The \newword{plactic monoid} $\Plactic(\cA)$ is the quotient of $\Free(\cA)$ by the \newword{Knuth relations}:
    \begin{align}
    acb &\sim cab \quad\text{for all } a\leq b < c;  \tag{K.1}\label{eq:P1}\\
      bca &\sim bac \quad \text{for all } a< b \leq c.  \tag{K.2}\label{eq:P2} 
    \end{align}
\end{definition}
The \newword{plactic Schur function} $\PlacticSchur_\nu$ is the image of $\FreeSchur_\nu$ under the projection map $\kappa : \Z\Free(\cA) \to \Z\Plactic(\cA)$. Classically, plactic Schur functions are usually studied with respect to the alphabet $\cN$ or one of its finite truncations; see \cite{Romik.Sniady:jdt} for consideration of the case $\cA = [0,1]_\mathbb{R}$ (as well as \cite{Romik.Sniady:bumping} for related discussion). We believe the case of arbitrary totally ordered alphabets deserves further attention.



\begin{remark}
    When the alphabet $\cA$ is infinite, the plactic Schur functions (resp.\ free Schur functions) are not elements of the monoid algebra of $\Plactic(\cA)$ (resp.\ $\Free(\cA)$) because they are infinite sums. Rather, one needs to consider a graded inverse limit as in Remark~\ref{rem:completion}. This is analogous to the situation with classical symmetric functions, which lie in the ring $\Z \llbracket x_1, x_2, \dots \rrbracket_{\rm gr}$ of power series of bounded degree.
\end{remark}

We now explain the motivation behind the definition of free Schur functions. This perspective will help us later to obtain analogous definitions in the shifted setting. 

\emph{Schensted insertion} \cite{Schensted} is an algorithm that produces a semistandard Young tableau $P(w)$ associated to a word $w \in \Free(\cA)$.
C.~Schensted \cite{Schensted} proved that 
    the number of columns in $P(w)$ is equal to the length of the longest weakly increasing subsequence of $w.$
C.~Greene \cite{Greene} established the stronger statement that the collection of increasing sequences of a word completely determines the shape of $P(w)$. This motivates a natural choice of representative for any tableau, first considered by Lascoux--Sch\"{u}tzenberger \cite{Lascoux.Schutzenberger:plaxique}. To wit, let $T \in \SSYT(\nu)$. The \newword{reading word} $\reading(T)$ of $T$ is obtained by reading the entries of $T$ by rows left to right and from bottom to top. Note that $\reading(T)$ is a representation of $T$ as a concatenation of weakly increasing segments. Moreover, the increasing decomposition $\reading(T)=w_1\ldots w_k$ is unique if $w_i$ is required to be of length $\nu_{\ell-i+1}$ and a longest increasing subword in $w_{i-1}w_i$ for all $i>1.$

Now, recalling that $s_\nu$ is a sum over $\SSYT(\nu)$, we can lift $s_\nu$ to $\Z\Free(\cA)$ by replacing the monomial ${\bf x}^{c(T)}$ corresponding to $T$ with $\reading(T) \in \Free(\cA)$. That is, we obtain the sum
\begin{equation}\label{eq:FreeSchurbyreading}
    \sum_{ T\in \SSYT(\nu)   }   \reading(T) \in \Z\Free(\cA).
\end{equation}
Noting that reading words of tableaux are exactly the words with the factorization properties considered above, we see that the element determined by \eqref{eq:FreeSchurbyreading} is exactly the free Schur function $\FreeSchur_\nu$ as defined in \eqref{eq:FreeSchur}.

\subsection{Shifted free Schur functions}
We now turn to giving analogous definitions in the shifted context.

    A \newword{hook subword} of a word $w$ in the alphabet $\cA$ is a subword $d\cdot i$ such that $d$ is a strictly decreasing subword and $i$ is a weakly increasing subword. If $w$ is a hook subword of itself, we say $w$ is a \newword{hook word}.
    Let $\nu$ be a strict partition of length $\ell$. Denote by $\hookword(\nu)$, the set of  words $w=w_1\ldots w_\ell$, where each $w_i$ is a hook word of length $\nu_{\ell-i+1}$ and such that, for all $i>1$, $w_i$ is a longest hook subword of $w_{i-1}w_i$.

   The following definition appears to be new, even for finite alphabets.
\begin{definition}\label{def:shifted_free}
        The \newword{shifted free Schur function} of shape $\nu$ is
\[\sFreeSchur_\nu \coloneqq \sum_{ w\: : \:  w\in \hookword(\nu)}   w.\]
\end{definition}  

\begin{example}\label{ex:21}
    \ytableausetup{boxsize=0.3em}
    Consider the shifted shape $\nu=(2,1)$. Then all concatenations of hook subwords of the lengths dictated by $\nu$ can be expressed in terms of variables $b<c<d$ standing for letters of $\cA$ as follows:

    If all letters are distinct, $bdc$ and $cdb$ is an exhaustive list of all words $w$ with the first two letters constituting a longest hook subword of $w$. Similarly, if exactly two letters are equal, $bdb$ and $ccb$ are all possible words $w\in \hookword(\nu)$. This analysis will be reused in the proof of Theorem~\ref{thm:main}.
\end{example}

We now turn to motivating Definition~\ref{def:shifted_free} and describing its connections to shifted tableau combinatorics.    
Any attempt to apply the philosophy of Section~\ref{sec:free!} to the task of defining  shifted (free) Schur functions needs to extend the notion of increasing segments to an appropriate analogue for shifted tableaux. Here, we are informed by the Sagan--Worley insertion algorithm, its relation to hook words, and the fact that the product defined from this algorithm does not define a plactic structure. 

To be more specific, since the product of two shifted tableaux $R$ and $T$ defined by inserting ({\it \`{a} la} Sagan--Worley) the reading word of $T$ into $R$ is not associative, one realizes that Sagan--Worley insertion does not define a shifted plactic structure. Its relation to hook words nonetheless makes evident that hook segments are the natural analogues of increasing segments in the shifted context, as made precise in the following theorem. 

\begin{theorem}[\cite{Worley,Sagan}]\label{thm:Sagan_Hook} 
    Let $T$ be the Sagan--Worley insertion tableau of a word $w$ in the alphabet $\cA$. Then the length of the first row of $T$ is equal to the size of the longest hook subword of $w$.
\end{theorem}


Note that although the connection between shifted tableaux and hook words is motivating for Definition~\ref{def:shifted_free}, it is not directly referenced in any way. Moreover, no consideration of Haiman's mixed insertion is involved in the definition.

\begin{remark}
    Despite being motivated in the opposite direction, starting from complete knowledge of mixed insertion and wishing to find a shifted plactic structure to realize it, the set of canonical representatives of a shifted shape $\nu$ was known to Serrano~\cite{Serrano} as ``shifted tableaux words.'' He also saw them as a natural generalization of the Lascoux--Sch\"{u}tzenberger representatives for unshifted tableaux.
\end{remark}



We now describe Haiman's \cite{Haiman} \emph{mixed insertion algorithm}.
Let $w=w_1\ldots w_n \in \Free(\cA)$ be a word without primed entries. Then, the \emph{mixed insertion} of $w$ is the shifted tableau $T$ constructed as follows.
Given a row $r$ of a tableau and given a value $z \in \cA$, we insert $z$ to $r$ by placing $z$ at the right end of $r$ if nothing in $r$ is strictly greater than $z$, while otherwise we find the smallest element $x$ in $r$ such that $x > z$ and insert $z$ in the position of $x$, while bumping out the value $x$. Similarly, given a column $c$ and a value $z' \in \cA' \setminus \cA$, we insert $z'$ to $c$ by placing $z$ at the bottom of $c$ if nothing in $c$ is strictly greater than $z'$, while otherwise we find the smallest element $x$ in $c$ such that $x > z'$ and insert $z'$ in the position of $x$, while bumping out the $x$.

\begin{itemize}
    \item Place $w_1$ in the first row of the diagram and set $i=2.$
    \item Insert $w_i$ into the first row as above and let $x$ be the entry that is bumped out (if any). 
    \item If $x$ is unprimed and not on the main diagonal, insert $x$ into the next row. If $x$ is primed or is unprimed but on the main diagonal, insert $x'$ into the next column. In either case, continue such row and column insertions on any further values that are bumped out.
    \item After the last insertion has been performed, set $i=i+1$ and repeat step $2$ if $i\leq n$; otherwise end and return the shifted tableau $T$.
\end{itemize}

\begin{figure}[htb]
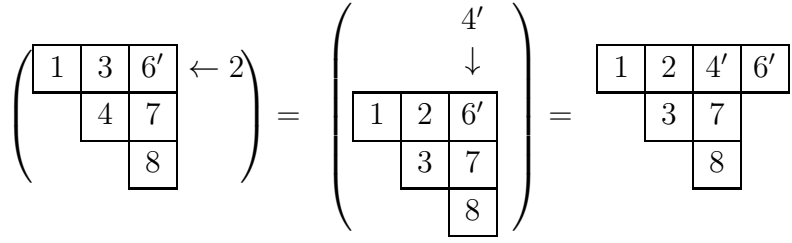

\ytableausetup{boxsize=normal,centertableaux}
    \centering
    \[
  \left( \begin{ytableau}
        1&3&6'&\none[\quad \leftarrow 2]\\
        \none&4&7 &\none\\
        \none & \none & 8 & \none
    \end{ytableau}  \hspace{2 mm} \right) =\hspace{2 mm}
    \left(
    \begin{ytableau}
        \none&\none&\none[4']\\
        \none&\none&\none[ \raisebox{3pt}{$\downarrow$}]\\
        1&2&6'\\
        \none&3&7 \\
        \none & \none & 8 
    \end{ytableau} \hspace{2 mm} \right) = \hspace{2 mm}
    \begin{ytableau}
        1&2&4'&6'\\
        \none&3&7 \\
        \none & \none & 8 
    \end{ytableau}
    \ytableausetup{nocentertableaux}
\]
    \caption{An example of the mixed insertion algorithm. The shown shifted tableaux illustrate the insertion of the element $2$ to the shifted semistandard Young tableau at left to obtain the rightmost tableau above.}
    \label{fig:shiftedlabel}
\end{figure}

Serrano \cite{Serrano} introduced the following monoid, which characterizes those words that have the same mixed insertion tableau. It is a remarkable fact that all of the necessary relations are in degree $4$. Theorem~\ref{thm:main} gives some new perspective on this fact; see also, Section~\ref{sec:staircase} for related discussion.

\begin{definition}[\cite{Serrano}]\label{def:shifted_plactic}
    The \newword{shifted plactic monoid} $\sPlactic(\cA)$ is the quotient of $\Free(\cA)$ by the \newword{shifted plactic relations}:
    \begin{align}
        &abdc\sim adbc \quad \text{for all }  a\leq b \leq c<d; \tag{SP.1}\label{eq:SP1} \\
    &acdb\sim acbd \quad \text{for all }  a\leq b < c \leq d; \tag{SP.2}\label{eq:SP2} \\ 
        &dacb\sim adcb \quad \text{for all }  a\leq b < c <d;\tag{SP.3}\label{eq:SP3} \\ 
        &badc\sim bdac \quad \text{for all }  a< b \leq c<d; \tag{SP.4}\label{eq:SP4} \\ 
        &cbda\sim cdba \quad \text{for all }  a< b <c \leq d; \tag{SP.5}\label{eq:SP5} \\ 
        &dbca\sim bdca \quad \text{for all }  a< b \leq c<d; \tag{SP.6}\label{eq:SP6} \\ 
        &bcda\sim bcad \quad \text{for all }  a< b \leq c\leq d; \label{eq:SP7} \tag{SP.7} \\ 
        &cadb\sim cdab \quad \text{for all }  a\leq  b < c\leq d. \tag{SP.8} \label{eq:SP8}
    \end{align}
\end{definition}

\section{Proof of Theorem~\ref{thm:LS}}\label{sec:LS_universal_proof}

\ytableausetup{boxsize=0.3em}

Fix an alphabet $\cA$. We will prove the theorem for $\LSCategory(\cA)$, often dropping $\cA$ from the notation for convenience.

Let $(\Monoid, \phi, \psi)$ be an arbitrary object in $\LSCategory$.
We now show that if $w, w' \in \Free$ are such that $\kappa(w) = \kappa(w')$ in $\Plactic$, then we also have $\phi(w)=\phi(w'),$ i.e.\ that the Knuth relations hold in $\Monoid$.

Consider the projections of the free Schur functions $\phi(\FreeSchur_{\ydiagram{1}})$ and $ \phi(\FreeSchur_{\ydiagram{1,1}})$. By condition \ref{plac2}, we have
\begin{align}
    \phi(\FreeSchur_{\ydiagram{1}})\phi(\FreeSchur_{\ydiagram{1,1}}) &= \phi(\FreeSchur_{\ydiagram{1,1}}) \phi(\FreeSchur_{\ydiagram{1}}) \label{eq:SchurCommute}
\end{align}
in $\mathbb{Z}\Monoid$. From \eqref{eq:SchurCommute}, we prove by cases that the Knuth relations hold in $\Monoid$. Note that this argument does not actually require one to know the Knuth relations in advance, although for conciseness we have used this foreknowledge here.

\medskip
\eqref{eq:P1}:
Let $a,b,c \in \cA$ with $a < b< c$, and consider the word $acb\in \Free$. For \eqref{eq:SchurCommute} to be satisfied in $\mathbb{Z}\Monoid$, it is necessary that $acb=a(cb)$ is equivalent to one of the words in Table~\ref{tab:11times1unshifted}, i.e., we must have $\phi(acb)=\phi(w)$ for one of the words $w$ in that table.
Moreover, by combining conditions \ref{plac1} and \ref{plac4}, $(\epsilon \circ \rho_I)(acb)=(\epsilon \circ \rho_I)(w)$ needs to hold for such a word $w$; it follows that $w$ has the same letters as $acb$, albeit in a possibly different order. 

Hence, $w$ is one of the words in Table~\ref{tab:abc_typeA}, and $acb$ must be equivalent to $(cb)a$, $(ca)b$, or $(ba)c$. (For reference, Table~\ref{tab:11times1unshifted} also shows the tableaux yielding these words; however, the tableaux are not necessary for the proof.) If either \[
\phi(acb)= \phi(cba) \quad \text{or} \quad \phi(acb)= \phi(bac),
\]
then, from condition \ref{plac4} on the interval $I=[a,b]$, we have correspondingly either
\begin{align*}
  \phi(ab)&=\phi \circ \rho_I(acb) = \phi \circ \rho_I(cba) = \phi(ba) \quad \text{or} \\
\phi(ab)&=\phi \circ \rho_I(acb)= \phi \circ \rho_I(bac) =\phi(ba).
\end{align*}
In either case, $\phi(ab) = \phi(ba)$.
However, if $\phi(ab)= \phi(ba)$, then by condition \ref{plac3}, we have that for all ordered morphisms $\omega$, $\phi\circ \omega (ab) = \phi \circ \omega (ba).$ That is, $\phi(xy)=\phi(yx)$ for all $x,y \in \cA$. 
This condition clearly implies that \eqref{eq:P1} holds in these cases.
If, on the other hand $\phi(acb)=\phi(cab)$, then that precisely says that the relation \eqref{eq:P1} holds in $\Monoid$, as well. 

Now, consider the case where $a = b < c$. Again, consider the word $aca \in \Free$. Since $aca$ is a monomial of $\FreeSchur_{\ydiagram{1}} \cdot \FreeSchur_{\ydiagram{1,1}}$, there must be a monomial $v$ in the product $\FreeSchur_{\ydiagram{1,1}} \cdot \FreeSchur_{\ydiagram{1}}$ with $\phi(aca) = \phi(v)$. Since then $\epsilon(aca) = \epsilon(v)$, $v$ must be a word of length $3$ using $a$ twice and $c$ once.
But Table~\ref{tab:aab_typeA} shows that there is a unique such monomial of $\FreeSchur_{\ydiagram{1,1}} \cdot \FreeSchur_{\ydiagram{1}}$. Thus, $v =caa$, in agreement with the $a=b$ case of \eqref{eq:P1}.

\medskip
\eqref{eq:P2}: 
Let $a,b,c\in \cA$ with $a<b<c$. Now, consider the word $b(ca)\in \Free.$ Since we have already shown that $\phi(acb) = \phi(cab)$, for \eqref{eq:SchurCommute} to hold in $\mathbb{Z}\Monoid$, it is necessary that either $\phi(bca)=\phi(cba)$ or $\phi(bca)=\phi(bac)$. 

If $\phi(bca)=\phi(cba)$, then by letting $I=[b,c]$, it follows that by condition \ref{plac4} that
\[\phi(bc) = \phi \circ \rho_I (bca) = \phi \circ \rho_I (cba) = \phi(cb) \]
so that again, by condition \ref{plac3}, we have that $\phi(bc)=\phi(cb)$ for all $b,c\in \cA.$ In particular, \eqref{eq:P2} is true in this case. 

If instead $\phi(bca)=\phi(bac)$, then this relation implies that \eqref{eq:P2} holds in $\Monoid$ too. 

Let us now consider the case where $a<b=c$. In this case, $bba$ is a monomial in the product $\sFreeSchur_{\ydiagram{1}}\cdot \sFreeSchur_{\ydiagram{1,1}}$, so that there must exist a monomial $v$ in the product $\sFreeSchur_{\ydiagram{1,1}}\cdot \sFreeSchur_{\ydiagram{1}}$ satisfying $\phi(bba)=\phi(v)$. 

Given that $v$ also needs to satisfy $\epsilon(bba)=\epsilon(v)$, $v$ contains two occurrences of $b$ and one of $a$, albeit possibly in a different order. However, Table~\ref{tab:abb_typeA} contains only one such monomial. Namely, $v$ is forced to be $bab$, agreeing with the $b=c$ case of \eqref{eq:P2}.

Thus, the Knuth relations hold in $\Monoid$.

It is straightforward to see that $(\Plactic,\kappa,\lambda) \in \LSCategory$.
Now, we can show that there is a unique morphism $\theta : \Plactic \to \Monoid$ for every $\Monoid \in \LSCategory$. Given an arbitrary element of $\Plactic$, since $\kappa$ is surjective, we may write the element as $\kappa(w)$ for some $w \in \Free$. Now, define $\theta(\kappa(w)) = \phi(w)$. This is well defined since if $\kappa(w) = \kappa(w')$, we have shown that $\phi(w) = \phi(w')$. Moreover, the left triangle of Equation~\eqref{eq:LS_rhombus} commutes by construction, and indeed $\theta$ is the only map that can make this triangle commute. That the right triangle commutes as well is straightforward from the fact that $\lambda$ and $\psi \circ \phi$ preserve the content of words. \qed



    \ytableausetup{boxsize=1.4em}
    
    \setlength{\arrayrulewidth}{0.01mm}
    \renewcommand{\arraystretch}{1.3}
     \ytableausetup{aligntableaux=center}
 \begin{table}[htbp]
    \begin{subtable}{0.29\textwidth}
		\begin{center}
		  \begin{tabular}{|c| c | c |}
			\hline
			Tableau       & Reading \\
			               & Word    \\
			\hline
			$
   \begin{ytableau}
			     b       \\
			     c       \\
			     \none \\
			    \end{ytableau}
       \bullet
   \begin{ytableau}
			a 
			\end{ytableau}

            $
			 & cba    \\
   
			\hline

			\begin{ytableau}
			a              \\
			c              \\
			   \none \\
			\end{ytableau} $\bullet$      
			\begin{ytableau}
			b 
			\end{ytableau}& cab    \\
			\hline

			\begin{ytableau}
			a               \\
			b               \\
			   \none \\
			\end{ytableau} $\bullet$ \begin{ytableau}
			c
			\end{ytableau} & bac    \\
			\hline
				    
		\end{tabular}
		    
		\end{center}
		\caption{All distinct entries.}
		\label{tab:abc_typeA}    
	\end{subtable} 
  \renewcommand{\arraystretch}{1.3}
	\begin{subtable}{0.29\textwidth}
		\begin{center}
		  \begin{tabular}{|c| c |}
			\hline
			Tableau & Reading \\
			         & Word    \\
			\hline
			\begin{ytableau}
			a           \\
			c           \\
			  \none \\
			\end{ytableau} $\bullet$
			\begin{ytableau}
			a 
			\end{ytableau}
			         & caa    \\
			\hline
						        
		\end{tabular}
		    
		\end{center}
		\caption{Two repeated entries with the smallest letter repeated.}
		\label{tab:aab_typeA}
  \end{subtable}
  \begin{subtable}{0.29\textwidth}
        \begin{center}
		      \begin{tabular}{|c|c|}
			\hline
			Tableau & Reading \\
			         & Word    \\
			\hline
			\begin{ytableau}
			  a      \\
			  b        \\
     \none \\
			\end{ytableau} $\bullet$
			\begin{ytableau}
			b 
			\end{ytableau}
			         & bab    \\
			\hline
						        
		\end{tabular}
		    
		\end{center}
		\caption{Two repeated entries with the biggest letter repeated.}
		\label{tab:abb_typeA}
	\end{subtable}
  \ytableausetup{boxsize=0.3em}
	\caption{The tableau products corresponding to $\FreeSchur_{\ydiagram{1,1}} \cdot \FreeSchur_{\ydiagram{1}}$, together with their reading words as considered in the proof of Theorem~\ref{thm:LS} in Section~\ref{sec:LS_universal_proof}.}\label{tab:11times1unshifted}
\end{table}

\section{Proof of Theorem~\ref{thm:main}}\label{sec:shifted_universal_proof}
Fix an alphabet $\cA$. We prove the theorem for $\Category(\cA)$, but routinely drop $\cA$ from the notation for conciseness. The main step is showing that every $\Monoid \in \Category$ satisfies the shifted Knuth relations, analogously to the proof of Theorem~\ref{thm:LS}, which relies on establishing the ordinary Knuth relations. As in that proof, this argument does not require foreknowledge of the shifted Knuth relations, but we refer to them here to streamline the analysis.

    Let $(\Monoid, \phi, \psi)$ be an arbitrary object of $\Category$. We now establish the shifted Knuth relations for $\Monoid$; that is, we show, for all $w,w' \in \Free$, that if $\sigma(w) = \sigma(w')$, then $\phi(w) = \phi(w')$.
    
    \ytableausetup{boxsize=0.3em}
    First, consider the shifted free Schur functions $\sFreeSchur_{\ydiagram{1}}, \sFreeSchur_{\ydiagram{2,1+1}} \in \Z \Free$. By \ref{splac2}, the images $\phi(\sFreeSchur_{\ydiagram{1}}), \phi(\sFreeSchur_{\ydiagram{2,1+1}}) \in \Z \Monoid$ commute. 
        We will consider each of the relations \eqref{eq:SP1}--\eqref{eq:SP8} of Definition~\ref{def:shifted_plactic} in turn. (As the arguments for the various relations are all very similar in format, we give progressively less detail as we proceed.) For the words contributing to the shifted free Schur function $\sFreeSchur_{\ydiagram{2,1+1}}$, see Example~\ref{ex:21}.



\medbreak
   \eqref{eq:SP1}:
   Let $a,b,c,d \in \cN$ with $a < b < c <d$.
   Consider $abdc = a\cdot bdc \in \Free$ Notice that $bdc$ is a monomial in the free Schur function $\sFreeSchur_{\ydiagram{2,1+1}}$. Hence, since there is no cancellation, $abdc$ is a monomial in the product $\sFreeSchur_{\ydiagram{1}}\sFreeSchur_{\ydiagram{2,1+1}}.$ Thus, 
   there must be a monomial $v$ in the product $\sFreeSchur_{\ydiagram{2,1+1}} \cdot \sFreeSchur_{\ydiagram{1}}$ such that
   $\phi(abdc) = \phi(v)$.
   Since $\epsilon(abdc) = \epsilon(v)$, $v$ must be a word of length exactly $4$, using the letters $a,b,c,d$ each exactly once. In the product $\sFreeSchur_{\ydiagram{2,1+1}}\cdot \sFreeSchur_{\ydiagram{1}} $, there are exactly $8$ monomials of degree $4$ using each of $a,b,c,d$ exactly once each; these are listed in Table~\ref{tab:abcstaircase}. (Again, for reference, Table~\ref{tab:allthetables} also shows the tableaux yielding these words; however, the tableaux are not necessary to the proof and the words can be constructed directly, as in Example~\ref{ex:21}.) Thus, $v$ must equal one of these $8$ words.
   

We will argue that, in fact, $v = adbc$. We consider the other $7$ possibilities from Table~\ref{tab:abcstaircase} in turn, and derive a contradiction for each. 
    
Suppose that $v= bdca$. Then $\phi(abdc) = \phi(bdca)$. By restricting to the interval $I = [a,b]$, we obtain that  $(\kappa \circ \rho_I) (abdc) = (\kappa \circ \rho_I)(bdca) \in \Plactic$, so that
\[
\kappa(ab) = \kappa(ba).
\]
This is a contradiction, as no two letters of the plactic monoid freely commute with each other.


    \ytableausetup{boxsize=1.4em}
    
    \setlength{\arrayrulewidth}{0.01mm}
    \renewcommand{\arraystretch}{1.3}
    
 \begin{table}[p]
    \begin{subtable}{0.34\textwidth}
		\begin{center}
		  \begin{tabular}{|c| c |}
			\hline
			Tableau       & Reading \\
			               & Word    \\
			\hline
			$
   \begin{ytableau}
			     b              & c       \\
			     \none          & d       \\
			     \none          & \none   \\
			    \end{ytableau}
       \bullet
   \begin{ytableau}
			a 
			\end{ytableau}

            $
			 & bdca    \\
   
			\hline

			\begin{ytableau}
			b              & c'      \\
			\none          & d       \\
			\none          & \none   
			\end{ytableau} $\bullet$      
			\begin{ytableau}
			a 
			\end{ytableau}& cdba    \\
			\hline

			\begin{ytableau}
			a              & c       \\
			\none          & d       \\
			\none          & \none   
			\end{ytableau} $\bullet$ \begin{ytableau}
			b 
			\end{ytableau} & adcb    \\
			\hline

			\begin{ytableau}
			a              & c'      \\
			\none          & d       \\
			\none          & \none   
			\end{ytableau} $\bullet$ \begin{ytableau}
			b 
			\end{ytableau} & cdab    \\
			\hline

			\begin{ytableau}
			a              & b       \\
			\none          & d       \\
			\none          & \none   
			\end{ytableau} $\bullet$ \begin{ytableau}
			c 
			\end{ytableau}
			 & adbc    \\
			\hline

			\begin{ytableau}
			a              & b'      \\
			\none          & d       \\
			\none          & \none   
			\end{ytableau} $\bullet$ \begin{ytableau}
			c 
			\end{ytableau} & bdac    \\
			\hline

			\begin{ytableau}
			a              & b       \\
			\none          & c       \\
			\none          & \none   
			\end{ytableau} $\bullet$ \begin{ytableau}
			d 
			\end{ytableau} & acbd    \\
			\hline

			\begin{ytableau}
			a              & b'      \\
			\none          & c       \\
			\none          & \none   
			\end{ytableau} $\bullet$ \begin{ytableau}
			d 
			\end{ytableau}& bcad    \\
			\hline
						    
		\end{tabular}
		    
		\end{center}
		\caption{All distinct entries.}
		\label{tab:abcstaircase}    
	\end{subtable} 
	\bigskip
	\hfill
  \renewcommand{\arraystretch}{1.3}
	\begin{subtable}{0.29\textwidth}
		\begin{center}
		  \begin{tabular}{|c| c |}
			\hline
			Tableau & Reading \\
			         & Word    \\
			\hline
			\begin{ytableau}
			b        & b       \\
			\none    & d       \\
			\none    & \none   
			\end{ytableau} $\bullet$
			\begin{ytableau}
			a 
			\end{ytableau}
			         & bdba    \\
			\hline
						        
			\begin{ytableau}
			a        & b       \\
			\none    & d       \\
			\none    & \none   
			\end{ytableau} $\bullet$
			\begin{ytableau}
			b 
			\end{ytableau}
			         & adbb    \\
			\hline
						        
			\begin{ytableau}
			a        & b'      \\
			\none    & d       \\
			\none    & \none   
			\end{ytableau} $\bullet$
			\begin{ytableau}
			b 
			\end{ytableau}
			         & bdab    \\
			\hline
						        
			\begin{ytableau}
			a        & b'      \\
			\none    & b       \\
			\none    & \none   
			\end{ytableau} $\bullet$
			\begin{ytableau}
			d 
			\end{ytableau}
			         & bbad    \\
			\hline
						        
		\end{tabular}
		    
		\end{center}
		\caption{Two repeated entries with the second smallest letter repeated.}
		\label{tab:abbstaircase}
		\vspace{0.2cm}
  
        \begin{center}
		      \begin{tabular}{|c|c|}
			\hline
			Tableau & Reading \\
			         & Word    \\
			\hline
			\begin{ytableau}
			b        & c'      \\
			\none    & c       \\
			\none    & \none   
			\end{ytableau} $\bullet$
			\begin{ytableau}
			a 
			\end{ytableau}
			         & ccba    \\
			\hline
						        
			\begin{ytableau}
			a        & c'      \\
			\none    & c       \\
			\none    & \none   
			\end{ytableau} $\bullet$
			\begin{ytableau}
			b 
			\end{ytableau}
			         & ccab    \\
			\hline
						        
			\begin{ytableau}
			a        & b       \\
			\none    & c       \\
			\none    & \none   
			\end{ytableau} $\bullet$
			\begin{ytableau}
			c 
			\end{ytableau}
			         & acbc    \\
			\hline
						        
			\begin{ytableau}
			a        & b'      \\
			\none    & c       \\
			\none    & \none   
			\end{ytableau} $\bullet$
			\begin{ytableau}
			c 
			\end{ytableau}
			         & bcac    \\
			\hline
						        
		\end{tabular}
		    
		\end{center}
		\caption{Two repeated entries with the biggest letter repeated.}
		\label{tab:addstaircase}
	\end{subtable}
	\bigskip
	\hfill
	\begin{subtable}{0.25\textwidth}
		
		\begin{tabular}{|c| c |}
			\hline
			Tableau & Reading \\
			         & Word    \\
			\hline
			\begin{ytableau}
			a        & c       \\
			\none    & d       \\
			\none    & \none   
			\end{ytableau} $\bullet$
			\begin{ytableau}
			a 
			\end{ytableau}
			         & adca    \\
			\hline
						        
			\begin{ytableau}
			a        & c'      \\
			\none    & d       \\
			\none    & \none   
			\end{ytableau} $\bullet$
			\begin{ytableau}
			a 
			\end{ytableau}
			         & cdaa    \\
			\hline
						        
			\begin{ytableau}
			a        & a       \\
			\none    & d       \\
			\none    & \none   
			\end{ytableau} $\bullet$
			\begin{ytableau}
			c 
			\end{ytableau}
			         & adac    \\
			\hline
						        
			\begin{ytableau}
			a        & a       \\
			\none    & c       \\
			\none    & \none   
			\end{ytableau} $\bullet$
			\begin{ytableau}
			d 
			\end{ytableau}
			         & acad    \\
			\hline
						        
		\end{tabular}
		\caption{Two repeated entries with the smallest letter repeated.}
		\label{tab:aabstaircase}
	\end{subtable} 
	\ytableausetup{boxsize=0.3em}
	\caption{The tableau products corresponding to $\sFreeSchur_{\ydiagram{2,1+1}} \cdot \sFreeSchur_{\ydiagram{1}}$, together with their mixed reading words as considered in the proof of Theorem~\ref{thm:main} in Section~\ref{sec:shifted_universal_proof}.}\label{tab:allthetables}
\end{table}

    For the same reason, 
    \[
    v \neq cdba, \quad v \neq bdac, \quad \text{and} \quad v \neq bcad,
    \]
    as each of these words has $b$ appearing before $a$.
    
    Using the same argument with the interval $[b,c]$ in place of the interval $[a,b]$, we conclude that 
     \[
    v \neq adcb, \quad v \neq cdab, \quad \text{and} \quad v \neq acbd,
    \]  as each of these words has $c$ appearing before $b$.
    This leaves 
    \[
    v = adbc
    \]
    as the only remaining possibility. The equivalence
    \[
     abdc \sim  adbc
    \]
    is
    shifted plactic equivalence \eqref{eq:SP1} of  Definition~\ref{def:shifted_plactic} for the case where all letters are distinct. 

    \ytableausetup{boxsize=0.3em}
    Now, consider the case where $a=b<c<d$. Hence, we consider the word $aadc = a \cdot adc \in \Free$. Note that this is a monomial appearing in $\sFreeSchur_{\ydiagram{1}}\cdot \sFreeSchur_{\ydiagram{2,1+1}}$. There must, therefore, be a monomial $v$ in the product $\sFreeSchur_{\ydiagram{2,1+1}} \cdot \sFreeSchur_{\ydiagram{1}}$ such that
   $\phi(aadc) = \phi(v)$.
   Since $\epsilon(aadc) = \epsilon(v)$, $v$ must be a word of length exactly $4$, using the letter $a$ exactly twice and the letters $c$ and $d$ each exactly once. In the product $\sFreeSchur_{\ydiagram{2,1+1}} \cdot \sFreeSchur_{\ydiagram{1}}$, there are exactly $4$ monomials of degree $4$ using each of $a,c,d$ the appropriate number of times. These monomials are listed in Table~\ref{tab:aabstaircase}. Thus, $v$ must equal one of these $4$ words.

    We show that $v = adac$ by eliminating the other three possibilities from Table~\ref{tab:aabstaircase}. First, suppose that $v = cdaa$. Then $\phi(aadc) = \phi(cdaa)$. By restricting to the interval $[c,d]$, we obtain that 
    $(\kappa \circ \rho_{[c,d]})(aadc) = (\kappa \circ \rho_{[c,d]})(cdaa) \in \Plactic$, whence $\kappa(dc) = \kappa(cd)$, contradicting again that no pair of letters in $\Plactic$ commute with each other. For the same reason, $v = acad$ is impossible.
    
   Finally, suppose that $v = adca$. Then $\phi(aada) = \phi(adca)$. By restricting to the interval $[a,c]$, we find that
   $(\kappa \circ \rho_{[a,c]})(aadc) = (\kappa \circ \rho_{[a,c]})(adca) \in \Plactic$, whence 
 $\kappa(aac) = \kappa(aca)$. This is a contradiction, as the relation $\kappa(xxy) = \kappa(xyx)$ only holds in $\Plactic$ for $x \geq y,$ whereas $a < c$.

    Having eliminated all other possibilities, it follows that $v = adac$, in agreement with the $a=b$ case of \eqref{eq:SP1}.
    
    We now consider the case where $a<b=c<d$. Hence, we consider the word $abdb = a \cdot bdb \in \Free$. This is a monomial in $\sFreeSchur_{\ydiagram{1}}\cdot \sFreeSchur_{\ydiagram{2,1+1}}$, and so there must be a monomial $v$ in the product 
    $\sFreeSchur_{\ydiagram{2,1+1}} \cdot \sFreeSchur_{\ydiagram{1}}$ such that
   $\phi(abdb) = \phi(v)$. Since $\epsilon(abdb) = \epsilon(v)$,
   $v$ must be a word of length $4$ using $b$ twice and $a,d$ each once. All four such monomials of $\sFreeSchur_{\ydiagram{2,1+1}} \cdot \sFreeSchur_{\ydiagram{1}}$ are listed in 
   Table~\ref{tab:abbstaircase}. 

   We show that $v = adbb$ by eliminating the other three possibilities from Table~\ref{tab:abbstaircase}. If $v = bdba$ or $v = bbad$, restricting to the interval $[a,b]$ shows that $\kappa(abb) = \kappa(bba)$, contradicting that no such relation holds in $\Plactic$ for $a < b$.
   If instead $v = bdab$, then restricting to the interval $[a,b]$ shows that $\kappa(abb) = \kappa(bab)$, which again is false in $\Plactic$. Thus, $v = adbb$, in agreement with the $b=c$ case of \eqref{eq:SP1}.
    
    Finally, we consider the case $a =b=c\leq d$. Hence, we consider the word $aada = a \cdot ada \in \Free$. This is a monomial in
    $\sFreeSchur_{\ydiagram{1}}\cdot \sFreeSchur_{\ydiagram{2,1+1}}$, and so there must be a monomial $v$ in the product 
    $\sFreeSchur_{\ydiagram{2,1+1}} \cdot \sFreeSchur_{\ydiagram{1}}$ such that
   $\phi(aada) = \phi(v)$. Since $\epsilon(aada) = \epsilon(v)$,
   $v$ must be a word of length $4$ using $a$ three times and $d$ once. The only such monomial of $\sFreeSchur_{\ydiagram{2,1+1}} \cdot \sFreeSchur_{\ydiagram{1}}$ is $adaa$. Thus, $v = adaa$, in agreement with the $a=b=c$ case of \eqref{eq:SP1}.
    
    \medbreak
   \eqref{eq:SP2}:
    First, we consider the case that $a <b <c <d$. Observe that $a\cdot cdb$ is a monomial in $\sFreeSchur_{\ydiagram{1}} \cdot \sFreeSchur_{\ydiagram{2,1+1}}$, so that it is equivalent to one of the terms $v$ in Table~\ref{tab:abcstaircase}. That is, there is a $v$ from Table~\ref{tab:abcstaircase} with $\phi(v)=\phi(acdb).$
    
    Restricting to the interval $[a,b]$ yields that
    \[
    v \neq bdca, \quad v \neq cdba, \quad v \neq bdac, \quad \text{and} \quad v\neq bcad.
    \]
    Restricting instead to $[b,c]$ yields $v\neq adbc$. Restricting to $[c,d]$ shows that $v \neq adcb$. 
    Suppose $v=cdab$. Then, the (trivial) restriction to the interval $I=[a,d]$ yields 
    \[
    \kappa(acdb) = \kappa(cdab) \in \Plactic,
    \] 
    which is incompatible with plactic equivalence, a contradiction.
    Thus, $v = acdb$, as desired and in accordance with the case $a<b<c<d$ of \eqref{eq:SP2}.

    Now, consider the case that $a=b<c < d$.
    Then we have that $a\cdot cda$ is a monomial of $\sFreeSchur_{\ydiagram{1}}\cdot \sFreeSchur_{\ydiagram{2,1+1}}$. Accordingly, it has to be equivalent under $\phi$ to a monomial $v$ in Table~\ref{tab:aabstaircase}.  
    
    Restricting to $[c,d]$ allows us to eliminate all possibilities besides either
    \[
        v = cdaa \quad \text{or} \quad v = acad.
    \]
    However, $v= cdaa$ is impossible because $\kappa(acda)\neq \kappa(cdaa) \in \Plactic$. Hence, $v = acad$, in accordance with the $a=b<c<d$ case of \eqref{eq:SP2}.

    Now, consider the case that $a<b<c = d$.
    The left side of \eqref{eq:SP2} is then $accb$, which must be equivalent under $\phi$ to one of the monomials $v$ in Table~\ref{tab:abbstaircase}.
    
    Restricting to $[a,b]$ eliminates the possibilities $v = ccb a$ and $v = bca c.$ Moreover, $v = ccab$ is impossible, since \[
\kappa(accb) \neq \kappa(ccab) \in \Plactic.
    \]
   Thus, $v = acb c$, in agreement with \eqref{eq:SP2} for the case $a<b<c = d$.
    
    Finally, consider the case $a=b<c = d$. Then, the left side of \eqref{eq:SP2} is $a cca$, and either
    \[v= cca  a \quad \text{or} \quad v= aca  c,\]
    since these are the only monomials $u$ in the product $\sFreeSchur_{\ydiagram{2,1+1}} \cdot \sFreeSchur_{\ydiagram{1}}$ with $\epsilon(u) = \epsilon(acca)$. 
    However, since 
    \[
    \kappa(acca) \neq \kappa(ccaa) \in \Plactic,
    \]
    we have
    $v= acac,$
    in accordance with the case $a=b<c = d$ of \eqref{eq:SP2}.

        \medbreak
   \eqref{eq:SP3}:
    Again, we begin with the case $a<b<c<d$. We identity what term $v$ of Table~\ref{tab:abcstaircase} satisfies $\phi(v)=\phi(d \cdot acb)$.
    Restricting to $[a,b]$ shows that 
    \[v= bdc\cdot a \quad \text{and} \quad v= cdb\cdot a \quad \text{and} \quad v= bda\cdot c \quad \text{and} \quad v= bca\cdot d\]
    are all impossible. 
    Furthermore, restricting to $[b,c]$ implies that $v= adb\cdot c$ is infeasible. Restricting to $[c,d]$ shows that $v=cda \cdot b$ and $v=acb\cdot d$ are also infeasible.
    The remaining equality $v=adcb$ is the desired equivalence in \eqref{eq:SP3} for the case of all letters distinct.

    Now, consider the case in which $a =b < c < d$.
    We obtain $d\cdot aca$ on the left side of \eqref{eq:SP3}. This must have the same image under $\phi$ as one of the words in Table~\ref{tab:aabstaircase}.
    Restricting to $[c,d]$ shows that  
    \[v \neq cdaa \quad \text{and} \quad v \neq aca  d.\]
    Moreover, $v \neq adac$, since $\kappa(daca) \neq \kappa(ada c) \in \Plactic$. Hence, $v = adca$, as desired, completing the case $a =b < c < d$ of \eqref{eq:SP3}.

        \medbreak
   \eqref{eq:SP4}:
    Suppose $a < b < c < d$. We need to identify an element $v$ in Table~\ref{tab:abcstaircase} with $\phi(v)=\phi(b\cdot adc)$.
    Restricting to $[a,b]$ eliminates the possibilities
    \[
    v=adc\cdot b, \quad v=cda\cdot b, \quad v=adb\cdot c, \quad \text{and} \quad v= acb\cdot d,
    \]
    while restricting to $[b,c]$ eliminates $v = cdb \cdot a$. Similarly, restricting to $[c,d]$ rules out $v= bca\cdot d$. Finally, $v\neq bdc\cdot a$, since $\kappa(badc)\neq \kappa(bdca)$. Hence, $v=bdac$, as desired, completing the case of \eqref{eq:SP4} with all distinct letters. 
    
    Now, consider the other case that $a < b=c < d$.
    Then, Table~\ref{tab:abbstaircase} provides the potential candidates for a word $v$ equivalent under $\phi$ to $b\cdot adb$. 
    
    Since $\phi(\sFreeSchur_{\ydiagram{2,1+1}} \cdot \sFreeSchur_{\ydiagram{1}}) = \phi(\sFreeSchur_{\ydiagram{1}}\cdot \sFreeSchur_{\ydiagram{2,1+1}})$, there must be a bijection between the monomials of $\sFreeSchur_{\ydiagram{2,1+1}} \cdot \sFreeSchur_{\ydiagram{1}}$ and $\sFreeSchur_{\ydiagram{1}}\cdot \sFreeSchur_{\ydiagram{2,1+1}}$ such that matched monomials $m_1, m_2$ satisfy $\phi(m_1) = \phi(m_2)$. We previously argued (in the $a<b=c<d$ case of \eqref{eq:SP1}) that $adb\cdot b$ was the only monomial $m$ in $\sFreeSchur_{\ydiagram{2,1+1}} \cdot \sFreeSchur_{\ydiagram{1}}$ with $\phi(abdb)=\phi(m)$. That is, it is the sole monomial that can be matched with $abdb$. Accordingly, $adbb$ cannot be the monomial $v$ matched in this case (or any case that follows) with $badb$.
    
    Also, by restriction to $[b,d]$, we find that $v \neq bbad$, since $\kappa(bdb)\neq \kappa(bbd) \in \Plactic$. Similarly, $v \neq bdb\cdot a$, since 
    \[
    \kappa(badb) \neq \kappa(bdba) \in \Plactic.
    \]
    Thus, $v = bdab$, as desired.

        \medbreak
   \eqref{eq:SP5}:
    Suppose $a < b < c < d$. We identity a word $v$ from Table~\ref{tab:abcstaircase} with $\phi(v)  = \phi(c\cdot bda)$.
    Restricting to $[a,b]$, we find that $v \neq adc\cdot b$, $v\neq cda\cdot b$, $v\neq acb\cdot d$, and $v\neq adb\cdot c$. Restricting to $[b,c]$, we further find that $v\neq bdc\cdot a$, $v \neq bda\cdot c$, and $v \neq bca\cdot d.$
    Thus, $v = cdb\cdot a$, as desired. 

    Now, consider the case $a < b < c = d$. The candidates for words $v$ equivalent under $\phi$ to $c\cdot bca$ are given in Table~\ref{tab:addstaircase}.
    Restriction to the interval $[a,b]$ shows that $v \neq cca\cdot b$ and $v \neq acb\cdot c$. 
    Finally, $v \neq bca\cdot c$, since 
    \[
    \kappa(cbca) \neq \kappa(bcac) \in \Plactic.
    \] Thus, $v = ccba$, as desired.

        \medbreak
   \eqref{eq:SP6}:
    Suppose $a < b < c < d$. 
    We identity a word $v$ from Table~\ref{tab:abcstaircase} with $\phi(v)  = \phi(d\cdot bca)$.
    There are three remaining possibilities from Table~\ref{tab:abcstaircase}, since the other five monomials have been matched in previous cases. Specifically, $adbc$ was used in \eqref{eq:SP1}; $acbd$ was matched in \eqref{eq:SP2}; $adcb$ was matched in \eqref{eq:SP3}; $bdac$ was used in \eqref{eq:SP4}; and $cdba$ was employed in \eqref{eq:SP5}.  
    By restriction to $[a,b]$, we find $v\neq cda\cdot b$. Similarly, by restriction to $[c,d]$, we find $v \neq bca\cdot d$. As all the other terms were employed for previous cases, we are left with $v= bdc\cdot a$, as desired. 
    
    Now, suppose $a < b = c < d$. Then, the word $d\cdot bba$ must be equivalent under $\phi$ to a word $v$ in Table~\ref{tab:abbstaircase}.
    Restriction to  $[a,b]$ shows that $v \neq adb\cdot b$. Restriction to $[b,d]$ shows that $v \neq bba\cdot d$. Finally, $v \neq bda\cdot b$, since 
    \[
    \kappa(dbba) \neq \kappa(bdab) \in \Plactic.
    \] Thus, $v = bdb\cdot a$, as desired. 

        \medbreak
   \eqref{eq:SP7}:
   Suppose $a < b < c < d$.
    We have two remaining monomials $v$ in Table~\ref{tab:abcstaircase} to match with $b\cdot cda$; since $bdca$ was used in \eqref{eq:SP6}, while $adbc$, $acbd$, $adcb$, $bdac$, and $cbda$ were used in earlier cases, as described in \eqref{eq:SP6}. 
    Of the remaining two monomials, the equality $v= cda\cdot b$ is impossible, as seen by restriction to $[a,b].$ Hence, $v= bca\cdot d$, as desired.
    
    Suppose $a < b = c < d$.
    Then, all but one of the words in Table~\ref{tab:abbstaircase} have been chosen; specifically, $adbb$ in \eqref{eq:SP1}, $bdab$ in \eqref{eq:SP4}, and $bdba$ in \eqref{eq:SP6}.
 The only remaining monomial yields the desired equality $v= bba\cdot d.$

    Suppose $a < b < c = d$.
    Then $b\cdot cca$ corresponds to one of the monomials $v$ in Table~\ref{tab:addstaircase}. 
    By restriction to $[a,b]$, the possibilities $v = cca\cdot b$ and $v = acb\cdot c$ can be ruled out. Moreover, by restriction to $[b,c]$, the equivalence $v = ccb\cdot a$ is also ruled out, since
    \[
    \kappa(bcc) \neq \kappa(ccb) \in \Plactic.
    \] Therefore, $v = bca\cdot c$, as desired. 

    Finally, suppose $a < b=c=d$. Then, $bbab$ is the only monomial $v$ in the sum $\sFreeSchur_{\ydiagram{2,1+1}} \cdot \sFreeSchur_{\ydiagram{1}}$ with $\epsilon(v) = \epsilon(b\cdot bba)$. Hence $v=bba\cdot b$, as desired.

        \medbreak
   \eqref{eq:SP8}:
    In this case, the desired equivalences are determined by the fact that there is only one remaining monomial available in each of the Tables~\ref{tab:abcstaircase},~\ref{tab:addstaircase}, and~\ref{tab:aabstaircase}. 
   
\medskip

It is straightforward to see from the definition that $(\sPlactic,\sigma,\pi) \in \Category$.
Now, we can show that there is a unique morphism $\theta : \sPlactic \to \Monoid$ for every $\Monoid \in \Category$. Given an arbitrary element of $\sPlactic$, since $\sigma$ is surjective, we may write the element as $\sigma(w)$ for some $w \in \Free$. Now, define $\theta(\sigma(w)) = \phi(w)$. This is well defined since if $\sigma(w) = \sigma(w')$, we have shown that $\phi(w) = \phi(w')$. Moreover, the left triangle of Equation~\eqref{eq:shifted_rhombus} commutes by construction, and indeed $\theta$ is the only map that can make this triangle commute. Commutativity of the right triangle follows from the fact that the shifted Knuth relations are implied by the ordinary Knuth relations, as may be seen by inspection of the lists of relations; that is, imposing the shifted Knuth relations and then the ordinary Knuth relations on $\Free$ is the same as imposing only the ordinary Knuth relations on $\Free$.
\qed

\section{Discussion on related axiomatizations}\label{sec:staircase}

Consider the axiom \ref{plac2}. It is natural to wonder what happens if $\phi(\FreeSchur_{\ydiagram{1,1}})$ is replaced by $\phi(\FreeSchur_{\ydiagram{2}})$, as this is the unique other partition of size $2$.  Tables \ref{tab:11times1unshifted} and \ref{tab:2times1unshifted} come to our aid, as they describe the terms of these free Schur functions. 

\ytableausetup{boxsize=1.4em}
    
    \setlength{\arrayrulewidth}{0.01mm}
    \renewcommand{\arraystretch}{1.3}
     \ytableausetup{aligntableaux=center}
 \begin{table}[htbp]
    \begin{subtable}{0.29\textwidth}
		\begin{center}
		  \begin{tabular}{|c| c | c |}
			\hline
			Tableau       & Reading \\
			               & Word    \\
			\hline
			$
   \begin{ytableau}
			     b   & c   
			    \end{ytableau}
       \bullet
   \begin{ytableau}
			a 
			\end{ytableau}

            $
			 & bca    \\
   
			\hline

			\begin{ytableau}
			a  &  c  
			\end{ytableau} $\bullet$      
			\begin{ytableau}
			b 
			\end{ytableau}& acb    \\
			\hline

			\begin{ytableau}
			a    & b
			\end{ytableau} $\bullet$ \begin{ytableau}
			c
			\end{ytableau} & abc    \\
			\hline
				    
		\end{tabular}
		    
		\end{center}
		\caption{All distinct entries.}
		\label{tab:abc_typeA_sec5}    
	\end{subtable} 
  \renewcommand{\arraystretch}{1.3}
	\begin{subtable}{0.29\textwidth}
		\begin{center}
		  \begin{tabular}{|c| c |}
			\hline
			Tableau & Reading \\
			         & Word    \\
			\hline
			\begin{ytableau}
			a & b 
			\end{ytableau} $\bullet$
			\begin{ytableau}
			a 
			\end{ytableau}
			         & aba    \\
			\hline
                \begin{ytableau}
			a & a 
			\end{ytableau} $\bullet$
			\begin{ytableau}
			b 
			\end{ytableau}
			         & aab    \\
			\hline
						        
		\end{tabular}
		    
		\end{center}
		\caption{Two repeated entries with the smallest letter repeated.}
		\label{tab:aab_typeA_sec5}
  \end{subtable}
  \begin{subtable}{0.29\textwidth}
        \begin{center}
		      \begin{tabular}{|c|c|}
			\hline
			Tableau & Reading \\
			         & Word    \\
			\hline
			\begin{ytableau}
			  a   & b
			\end{ytableau} $\bullet$
			\begin{ytableau}
			b 
			\end{ytableau}
			         & abb    \\
			\hline
			\begin{ytableau}
			  b   & b
			\end{ytableau} $\bullet$
			\begin{ytableau}
			a 
			\end{ytableau}
			         & bba    \\
			\hline			        
		\end{tabular}
		    
		\end{center}
		\caption{Two repeated entries with the biggest letter repeated.}
		\label{tab:abb_typeA_sec5}
	\end{subtable}
  \ytableausetup{boxsize=0.3em}
	\caption{The tableau products corresponding to $\FreeSchur_{\ydiagram{2}} \cdot \FreeSchur_{\ydiagram{1}}$
, as considered in Section~\ref{sec:staircase}.}\label{tab:2times1unshifted}
\end{table}

A comparison of these tables reveals that after cancelling equal terms on both sides, $\phi(\FreeSchur_{\ydiagram{1}})\cdot \phi(\FreeSchur_{\ydiagram{1,1}}) = \phi(\FreeSchur_{\ydiagram{1,1}}) \cdot \phi(\FreeSchur_{\ydiagram{1}})$ imposes the same conditions as $\phi(\FreeSchur_{\ydiagram{1}})\cdot \phi(\FreeSchur_{\ydiagram{2}})=\phi(\FreeSchur_{\ydiagram{2}}) \cdot \phi(\FreeSchur_{\ydiagram{1}}).$ Thus we have the following.

\begin{proposition}
    The category defined by replacing $\phi(\FreeSchur_{\ydiagram{1,1}})$ by $\phi(\FreeSchur_{\ydiagram{2}})$ in \ref{plac2} is equal to $\LSCategory(\cA).$ In particular, $\Plactic(\cA)$ is the initial object of both. \qed
\end{proposition}

Similarly, it is interesting to ask what happens to $\Category$ if \ref{splac2} is modified in an analogous way. This time, we have to consider two different variations: Replacing $\phi(\sFreeSchur_{\ydiagram{2,1+1}})$ either by 
$\phi(\sFreeSchur_{\ydiagram{2}})$ or by $\phi(\sFreeSchur_{\ydiagram{3}})$. Note also, that these are the only other partitions of sizes $2$ and $3$.

In the first case, it suffices to observe that any word of length $2$ is a hook word. Hence, $\sFreeSchur_{\ydiagram{2}}$ is a sum of all possible words of length $2$. Accordingly, $\sFreeSchur_{\ydiagram{1}}\cdot \sFreeSchur_{\ydiagram{2}} = \sFreeSchur_{\ydiagram{2}} \cdot \sFreeSchur_{\ydiagram{1}} \in \Z\Free(\cA)$. Thus, modifying \ref{splac2} in this way makes the axiom redundant. In particular, the initial object of this modified category is the free monoid $\Free(\cA).$

For the second case, we compare Table~\ref{tab:allthetables} to the Tables~\ref{tab:all_the_3by1_tables} and~\ref{tab:bcc}. As in the unshifted case, after cancelling identical words on both sides of $\phi(\sFreeSchur_{\ydiagram{1}})\cdot\phi(\sFreeSchur_{\ydiagram{3}}) = \phi(\sFreeSchur_{\ydiagram{3}}) \cdot \phi(\sFreeSchur_{\ydiagram{1}})$ and $\phi(\sFreeSchur_{\ydiagram{1}}) \cdot\phi(\sFreeSchur_{\ydiagram{2,1+1}}) = \phi(\sFreeSchur_{\ydiagram{2,1+1}}) \cdot \phi(\sFreeSchur_{\ydiagram{1}})$ we see that both commutations impose the same relations. Thus, we have the following proposition.

\begin{proposition}
    The category defined by replacing $\phi(\FreeSchur_{\ydiagram{2,1+1}})$ by $\phi(\FreeSchur_{\ydiagram{3}})$ in \ref{splac2} is equal to $\Category(\cA).$ In particular, $\sPlactic(\cA)$ is the initial object of both. \qed
\end{proposition}

 \ytableausetup{centertableaux}

\ytableausetup{boxsize=normal}
    \setlength{\arrayrulewidth}{0.01mm}
  
    \begin{table}[h]
    \begin{center}
    \renewcommand{\arraystretch}{2}
    \begin{subtable}[h!]{0.25\textwidth}
        \begin{center}
        \begin{tabular}{|c| c |}
        \hline
            $P_{(1)}P_{(3)}$ & $P_{(3)}P_{(1)}$ \\
            \hline
            abcd & bcda \\
        \hline
        
          acbd & cbda \\
        \hline
          adbc & dbca \\
        \hline
          adcb & dcba \\
    \hline
    
         bacd & acdb \\
        \hline
    bcad & cadb \\
        \hline
    bdac & dacb \\
        \hline
    bdca & dcab \\
        \hline
        cabd & abdc \\
        \hline
        cbad & badc \\
        \hline
        cdab & dabc \\
        \hline
        cdba & dbac \\
        \hline
        dabc & abcd \\
        \hline
        dbac & bacd \\
        \hline
        dcab & cabd \\
        \hline
        dcba & cbad \\
        \hline
        \end{tabular}
        \caption{All distinct entries.}
        \label{tab:aab}
        \end{center}
    \end{subtable}
    \hfill
    \renewcommand{\arraystretch}{2.5}
    \begin{subtable}[h!]{0.4\textwidth}
        \begin{center}
        \begin{tabular}{|c| c |}
        \hline
            Tableau & Reading Word \\
            \hline
            \begin{ytableau}
            a & b' & c'
    \end{ytableau} $\bullet$
            \begin{ytableau}
                   a 
            \end{ytableau}
        & cbaa \\
        \hline
        
            \begin{ytableau}
            a & b'& c
    \end{ytableau} $\bullet$
            \begin{ytableau}
                   a 
            \end{ytableau}
        & baca \\
        \hline
            \begin{ytableau}
            a & b &c' 
    \end{ytableau} $\bullet$
            \begin{ytableau}
                   a 
            \end{ytableau}
        & caba \\
        \hline
        
            \begin{ytableau}
            a & b & c 
    \end{ytableau} $\bullet$
            \begin{ytableau}
                   a 
            \end{ytableau}
        & abca \\
        \hline

            \begin{ytableau}
            a & a & c 
    \end{ytableau} $\bullet$
            \begin{ytableau}
                   b 
            \end{ytableau}
        & aacb \\
        \hline

            \begin{ytableau}
            a & a & c' 
    \end{ytableau} $\bullet$
            \begin{ytableau}
                   b 
            \end{ytableau}
        & caab \\
        \hline

            \begin{ytableau}
            a & a & b 
    \end{ytableau} $\bullet$
            \begin{ytableau}
                   c 
            \end{ytableau}
        & aabc \\
        \hline

            \begin{ytableau}
            a & a & b' 
    \end{ytableau} $\bullet$
            \begin{ytableau}
                   c 
            \end{ytableau}
        & baac \\
        \hline
        
        \end{tabular}
        \caption{Two repeated entries with the smallest letter repeated.}
        \label{tab:abb}
        \end{center}
    \end{subtable}
    \hfill
    \renewcommand{\arraystretch}{2.5}
    \begin{subtable}[h!]{0.27\textwidth}
        \begin{center}
        \begin{tabular}{|c| c |}
        \hline
            Tableau & Reading Word \\
            \hline
            \begin{ytableau}
            b & b & c'
    \end{ytableau} $\bullet$
            \begin{ytableau}
                   a 
            \end{ytableau}
        & bbca \\
        \hline
        
            \begin{ytableau}
            b & b& c' 
    \end{ytableau} $\bullet$
            \begin{ytableau}
                   a 
            \end{ytableau}
        & cbba \\
        \hline
        
            \begin{ytableau}
            a & b &c 
    \end{ytableau} $\bullet$
            \begin{ytableau}
                   b 
            \end{ytableau}
        & abcb \\
        \hline
        
            \begin{ytableau}
            a & b' & c 
    \end{ytableau} $\bullet$
            \begin{ytableau}
                   b 
            \end{ytableau}
        & bacb \\
        \hline

            \begin{ytableau}
            a & b & c' 
    \end{ytableau} $\bullet$
            \begin{ytableau}
                   b 
            \end{ytableau}
        & cabb \\
        \hline

            \begin{ytableau}
            a & b' & c' 
    \end{ytableau} $\bullet$
            \begin{ytableau}
                   b 
            \end{ytableau}
        & cbab \\
        \hline

            \begin{ytableau}
            a & b & b 
    \end{ytableau} $\bullet$
            \begin{ytableau}
                   c 
            \end{ytableau}
        & abbc \\
        \hline

            \begin{ytableau}
            a & b' & b 
    \end{ytableau} $\bullet$
            \begin{ytableau}
                   c 
            \end{ytableau}
        & babc \\
        \hline
        
        \end{tabular}
        \caption{Two repeated entries with the second smallest letter repeated.}
        \label{tab:abd}
        \end{center}
    \end{subtable}
    \ytableausetup{boxsize=0.3em}
    \caption{The tableau products corresponding to $\sFreeSchur_{\ydiagram{3}} \cdot \sFreeSchur_{\ydiagram{1}}$, together with their reading words as considered in Section~\ref{sec:staircase}. See Figure~\ref{tab:bcc} for the remaining case.}\label{tab:all_the_3by1_tables}
    \end{center}
    \end{table}

\begin{table}[h!]
\ytableausetup{boxsize=1.4em}
    \renewcommand{\arraystretch}{2}
        \begin{center}
        \begin{tabular}{|c| c |}
        \hline
            Tableau & Reading Word \\
            \hline
            \begin{ytableau}
            b & c & c
    \end{ytableau} $\bullet$
            \begin{ytableau}
                   a 
            \end{ytableau}
        & bcaa \\
        \hline
        
            \begin{ytableau}
            b & c'& c 
    \end{ytableau} $\bullet$
            \begin{ytableau}
                   a 
            \end{ytableau}
        & cbca \\
        \hline
        
            \begin{ytableau}
            a & c &c 
    \end{ytableau} $\bullet$
            \begin{ytableau}
                   b 
            \end{ytableau}
        & accb \\
        \hline
        
            \begin{ytableau}
            a & c' & c 
    \end{ytableau} $\bullet$
            \begin{ytableau}
                   b 
            \end{ytableau}
        & cacb \\
        \hline

            \begin{ytableau}
            a & b & c 
    \end{ytableau} $\bullet$
            \begin{ytableau}
                   c 
            \end{ytableau}
        & abcc \\
        \hline

            \begin{ytableau}
            a & b' & c 
    \end{ytableau} $\bullet$
            \begin{ytableau}
                   c 
            \end{ytableau}
        & bacc \\
        \hline

            \begin{ytableau}
            a & b & c' 
    \end{ytableau} $\bullet$
            \begin{ytableau}
                   c 
            \end{ytableau}
        & cabc \\
        \hline

            \begin{ytableau}
            a & b' & c' 
    \end{ytableau} $\bullet$
            \begin{ytableau}
                   c 
            \end{ytableau}
        & cbac \\
        \hline
        
        \end{tabular}
       \ytableausetup{boxsize=0.3em}
        \caption{The remaining case of tableau products corresponding to $\sFreeSchur_{\ydiagram{3}} \cdot \sFreeSchur_{\ydiagram{1}}$ from Table~\ref{tab:all_the_3by1_tables}. This case is where there are two repeated entries with the biggest letter repeated.}
        \label{tab:bcc}
        \end{center}
    \end{table}

\section*{Acknowledgements}
The authors are grateful for useful conversations with Alejandro Morales and Luis Serrano.

Both authors were partially supported by a Discovery Grant (RGPIN-2021-02391) and Launch Supplement (DGECR-2021-00010) from
the Natural Sciences and Engineering Research Council of Canada. SE also acknowledges partial support from a Sinclair Graduate Scholarship from the University of Waterloo.

\bibliographystyle{amsalpha} 
\bibliography{shifted.bib}

\providecommand{\bysame}{\leavevmode\hbox to3em{\hrulefill}\thinspace}
\providecommand{\MR}{\relax\ifhmode\unskip\space\fi MR }
\providecommand{\MRhref}[2]{%
  \href{http://www.ams.org/mathscinet-getitem?mr=#1}{#2}
}
\providecommand{\href}[2]{#2}
\begin{thebibliography}{Wor84}

\bibitem[FG98]{Fomin.Greene}
Sergey Fomin and Curtis Greene, \emph{Noncommutative {S}chur functions and
  their applications}, Discrete Math. \textbf{193} (1998), no.~1-3, 179--200,
  Selected papers in honor of Adriano Garsia (Taormina, 1994). \MR{1661368}

\bibitem[Gil19]{Gillespie}
Maria Gillespie, \emph{Variations on a theme of {S}chubert calculus}, Recent
  trends in algebraic combinatorics, Assoc. Women Math. Ser., vol.~16,
  Springer, Cham, 2019, pp.~115--158. \MR{3969573}

\bibitem[Gre74]{Greene}
Curtis Greene, \emph{An extension of {S}chensted's theorem}, Advances in Math.
  \textbf{14} (1974), 254--265.

\bibitem[Hai89]{Haiman}
Mark~D. Haiman, \emph{On mixed insertion, symmetry, and shifted {Y}oung
  tableaux}, J. Combin. Theory Ser. A \textbf{50} (1989), no.~2, 196--225.

\bibitem[HH92]{Hoffman.Humphreys}
P.~N. Hoffman and J.~F. Humphreys, \emph{Projective representations of the
  symmetric groups: ${Q}$-functions and shifted tableaux}, Oxford Mathematical
  Monographs, Oxford Science Publications, The Clarendon Press, Oxford
  University Press, New York, 1992.

\bibitem[Knu70]{Knuth}
Donald~E. Knuth, \emph{Permutations, matrices, and generalized {Y}oung
  tableaux}, Pacific J. Math. \textbf{34} (1970), 709--727.

\bibitem[Lot02]{Lothaire}
M.~Lothaire, \emph{Algebraic combinatorics on words}, Encyclopedia of
  Mathematics and its Applications, vol.~90, Cambridge University Press,
  Cambridge, 2002, A collective work by Jean Berstel, Dominique Perrin, Patrice
  Seebold, Julien Cassaigne, Aldo De Luca, Steffano Varricchio, Alain Lascoux,
  Bernard Leclerc, Jean-Yves Thibon, Veronique Bruyere, Christiane Frougny,
  Filippo Mignosi, Antonio Restivo, Christophe Reutenauer, Dominique Foata,
  Guo-Niu Han, Jacques Desarmenien, Volker Diekert, Tero Harju, Juhani
  Karhumaki and Wojciech Plandowski, With a preface by Berstel and Perrin.

\bibitem[LS81]{Lascoux.Schutzenberger:plaxique}
Alain Lascoux and Marcel-Paul Sch\"{u}tzenberger, \emph{Le mono\"{\i}de
  plaxique}, Noncommutative structures in algebra and geometric combinatorics
  ({N}aples, 1978), Quad. ``Ricerca Sci.'', vol. 109, CNR, Rome, 1981,
  pp.~129--156.

\bibitem[Nov00]{Novelli}
Jean-Christophe Novelli, \emph{On the hypoplactic monoid}, Discrete Math.
  \textbf{217} (2000), no.~1-3, 315--336, Formal power series and algebraic
  combinatorics (Vienna, 1997).

\bibitem[RS15]{Romik.Sniady:jdt}
Dan Romik and Piotr \'{S}niady, \emph{Jeu de taquin dynamics on infinite
  {Y}oung tableaux and second class particles}, Ann. Probab. \textbf{43}
  (2015), no.~2, 682--737. \MR{3306003}

\bibitem[RS16]{Romik.Sniady:bumping}
\bysame, \emph{Limit shapes of bumping routes in the {R}obinson-{S}chensted
  correspondence}, Random Structures Algorithms \textbf{48} (2016), no.~1,
  171--182. \MR{3432576}

\bibitem[Sag87]{Sagan}
Bruce~E. Sagan, \emph{Shifted tableaux, {S}chur {$Q$}-functions, and a
  conjecture of {R}. {S}tanley}, J. Combin. Theory Ser. A \textbf{45} (1987),
  no.~1, 62--103.

\bibitem[Sch61]{Schensted}
C.~Schensted, \emph{Longest increasing and decreasing subsequences}, Canadian
  J. Math. \textbf{13} (1961), 179--191. \MR{121305}

\bibitem[Sch77]{Schutzenberger}
M.-P. Sch\"{u}tzenberger, \emph{La correspondance de {R}obinson}, Combinatoire
  et repr\'{e}sentation du groupe sym\'{e}trique ({A}ctes {T}able {R}onde
  {CNRS}, {U}niv. {L}ouis-{P}asteur {S}trasbourg, {S}trasbourg, 1976), Lecture
  Notes in Math., Vol. 579, Springer, Berlin-New York, 1977, pp.~59--113.

\bibitem[Sch97]{Schutzenberger:pourleplaxique}
Marcel-Paul Sch\"{u}tzenberger, \emph{Pour le mono\"{\i}de plaxique}, Math.
  Inform. Sci. Humaines (1997), no.~140, 5--10.

\bibitem[Ser10]{Serrano}
Luis Serrano, \emph{The shifted plactic monoid}, Math. Z. \textbf{266} (2010),
  no.~2, 363--392.

\bibitem[Wor84]{Worley}
Dale~Raymond Worley, \emph{A theory of shifted {Y}oung tableaux}, ProQuest LLC,
  Ann Arbor, MI, 1984, Thesis (Ph.D.)--Massachusetts Institute of Technology.

\end{thebibliography}
\end{document}